\documentclass[leqno,12pt]{amsart}
\usepackage{amsmath,amstext,amssymb,amsopn,amsthm,mathrsfs}

\theoremstyle{plain}
\newtheorem{theorem}{Theorem}[section]
\newtheorem{defi}{Definition}[section]

\newtheorem{lemma}{Lemma}[section]
\newtheorem{proposition}{Proposition}[section]
\newtheorem{obs}{Observation}[section]
\numberwithin{equation}{section}
\newcommand{\dem}{\medskip \par \noindent \mbox{\bf Proof. }}
\def\ep{\hfill{$\Box $}}

\begin{document}

\title[On Gaussian Lipschitz spaces ]{On Gaussian Lipschitz spaces and the boundedness of Fractional Integrals and Fractional Derivatives  on them.}

\author{A. Eduardo Gatto}
\address{Department of Mathematical Sciences, DePaul University, Chicago, IL,
   60614, USA.}
   \email{aegatto@depaul.edu}
\author{Wilfredo~O.~Urbina~R.}
\address{Department of Mathematics and Actuarial Sciences, Roosevelt University, 430 S. Michigan Ave. Chicago, IL, 60605, USA}
\email{wurbinaromero@roosevelt.edu}



\begin{abstract}
The main purpose of this paper is to study the boundedness of Gaussian fractional integrals and derivatives associated to Hermite polynomial expansions on  Gaussian Lipschitz spaces $Lip_{\alpha}(\gamma)$. To get these results we introduce formulas for these operators in terms of the Hermite-Poisson semigroup as well as the Gaussian Lipschitz spaces. This approach was originally developed for the classical  Poisson integral. These proofs  can also be extended to the case of Laguerre and Jacobi expansions. In subsequent papers we will study the same operators on Gaussian Besov-Lipschitz and Triebel-Lizorkin spaces.

\end{abstract}

\maketitle

\section{Preliminaries.}
Let us consider $\mathbb{R}^d$, with the Gaussian measure $d\gamma(x)=\frac{e^{-\left|x\right|^2}
}{\pi^{d/2}}\,dx$, $x\in\mathbb{R}^d$ and the Ornstein-Uhlenbeck
differential operator
\begin{equation}\label{OUop}
L=\frac12\triangle_x-\left\langle x,\nabla _x\right\rangle.
\end{equation}

Let $\nu=(\nu _1,...,\nu_d)$ be a
multi-index, where $\nu_i$ is a non-negative integer, $1\leq i \leq d$; let $\nu
!=\prod_{i=1}^d\nu _i!,$ $\left| \nu \right| =\sum_{i=1}^d\nu _i,$ $%
\partial _i=\frac \partial {\partial x_i},$ for each $1\leq i\leq d$ and $%
\partial ^\nu =\partial _1^{\nu _1}...\partial _d^{\nu_d}.$

Let  $h_\nu $ be the Hermite polynomial of order
$\nu$, in $d$ variables,
\begin{equation}
h_\nu (x)=\frac 1{\left( 2^{\left| \nu \right| }\nu
!\right)
^{1/2}}\prod_{i=1}^d(-1)^{\nu _i}e^{x_i^2}\frac{\partial ^{\nu _i}}{%
\partial x_i^{\nu _i}}(e^{-x_i^2}),
\end{equation}
then, it is well known, that the Hermite polynomials are
eigenfunctions of $L$,
\begin{equation}\label{eigen}
L h_{\nu}(x)=-\left|\nu \right|h_\nu(x).
\end{equation}
Given a function $f$ $\in L^1(\gamma )$ its
$\nu$-Fourier-Hermite coefficient is defined by
\begin{equation}
\hat{f}(\nu) =<f, h_\nu>_{\gamma}
=\int_{\mathbb{R}^d}f(x)h_\nu (x) \, d\gamma(x).
\end{equation}

For each $n$, a non negative integer, let $C_n$ be the closed subspace of $L^2(\gamma)$ generated by
the linear combinations of $\left\{ h_\nu \ :\left| \nu
\right| =n\right\}$. By the orthogonality of the Hermite
polynomials with respect to $\gamma$, it is easy to see that
$\{C_n\}$ is an orthogonal decomposition of $L^2(\gamma)$,
$$ L^2(\gamma) = \bigoplus_{n=0}^{\infty} C_n$$
which is called the Wiener chaos.

If $J_n$ is the orthogonal projection  of $L^2(\gamma)$ onto
$C_n$, that is, for $f\in L^2(\gamma)$
\[
J_n f=\sum_{\left|\nu \right|=n}\hat{f}(\nu) h_\nu,
\]
then the Hermite expansion of $f$ can be written as
$$ f = \sum_{n=0}^{\infty} J_n f.$$
\\
Let us consider the Ornstein-Uhlenbeck semigroup $\left\{
T_t\right\} _{t\geq 0}$, i.e.
\begin{eqnarray}\label{01}
\nonumber T_t f(x)&=&\frac 1{\left( 1-e^{-2t}\right) ^{d/2}}\int_{\mathbb{R}^d}e^{-\frac{%
e^{-2t}(\left| x\right| ^2+\left| y\right| ^2)-2e^{-t}\left\langle
x,y\right\rangle }{1-e^{-2t}}}f(y)\, d\gamma(y)\\
& = & \frac{1}{\pi^{d/2}(1-e^{-2t})^{d/2}}\int_{\mathbb R^d} e^{-
\frac{|y-e^{-t}x|^2}{1-e^{-2t}}} f(y) dy,
\end{eqnarray}
for $t  > 0$ and $T_0  f= f.$\\
The family $\left\{ T_t\right\}_{t\geq 0}$ is a strongly
continuous positive symmetric diffusion semigroup, see \cite{st2} on $L^p(\gamma)$, $1 \leq p <
\infty$, with infinitesimal generator $L$. Now, by Bochner subordination formula, see Stein \cite{se70},  the Poisson-Hermite 
semigroup $\left\{ P_t\right\} _{t\geq0}$ is defined  as
\begin{equation}\label{PoissonH}
\nonumber P_t f(x)=\frac 1{\sqrt{\pi }}\int_0^{\infty} \frac{e^{-u}}{\sqrt{u}}T_{t^2/4u}f(x)du.
\end{equation}
In addition, changing the variable, we have the following representation
$$P_t f(x)=\int_0^{\infty} T_s f(x) \mu^{(1/2)}_t(ds),$$
where the measure
\begin{equation}\label{onesided1/2}
\mu^{(1/2)}_t(ds) = \frac t{2\sqrt{\pi
}}\frac{e^{-t^2/4s}}{s^{3/2}}ds = g(t,s) ds,
\end{equation}
is called the one-side stable measure on $(0, \infty)$ of order
$1/2$.

 From (\ref{01}) we
obtain, after the change of variable $r=e^{-t^2/4u}$,
\begin{eqnarray}\label{03}
\nonumber P_t f(x)&=&\frac 1{2\pi
^{(d+1)/2}}\int_{\mathbb{R}^d}\int_0^1t\frac{\exp \left( t^2/4\log
r\right) }{(-\log r)^{3/2}}\frac{\exp \left( \frac{-\left|
y-rx\right| ^2}{1-r^2}\right) }{(1-r^2)^{d/2}}\frac{dr}rf(y)dy\\
&=& \int_{\mathbb{R}^d} p(t,x,y) f(y)dy,
\end{eqnarray}
with
\begin{equation}
p(t,x,y) = \frac 1{2\pi ^{(d+1)/2}}\int_0^1t\frac{\exp \left(
t^2/4\log r\right) }{(-\log r)^{3/2}}\frac{\exp \left(
\frac{-\left| y-rx\right| ^2}{1-r^2}\right)
}{(1-r^2)^{d/2}}\frac{dr}r.
\end{equation}

The family $\left\{ P_t\right\}_{t \geq 0}$ is also a strongly
continuous positive symmetric diffusion semigroup on $L^p(\gamma)$, $1 \leq p < \infty$,
with infinitesimal generator $-(-L)^{1/2}$. In particular $P_t 1 =1$.

Observe that by (\ref{eigen}), $h_{\nu}$ are eigenfunctions of  $\left\{ T_t\right\}_{t \geq 0}$ and $\left\{ P_t\right\}_{t \geq 0}$we have that
\begin{equation}\label{OUHerm}
T_t h_\nu(x)=e^{-t\left| \nu\right|}h_\nu(x),
\end{equation}
and
\begin{equation}\label{PHHerm}
 P_t h_\nu(x)=e^{-t\sqrt{\left| \nu\right|}}h_\nu(x)
\end{equation}
and therefore, if $f \in L^2(\gamma)$ 
$$ T_t  f =  \sum_n e^{-nt} J_n f \quad
\mbox{and} \quad  P_t  f =  \sum_n e^{-\sqrt{n}\,t} J_n f.$$
Since $\left\{P_t\right\}_{t\geq 0}$, is a strongly continuous semigroup we have
\begin{equation}\label{d1}
\lim_{t \to 0^{+}}P_t f(x)=f(x), \quad \mbox{a.e.} \; x,
\end{equation}
and we have,
\begin{equation} \label{d2}
\lim_{t \to \infty} P_t f(x) = \int_{\mathbb{R}^d} f(y) d\gamma (y), \quad \mbox{a.e.} \; x. \\
\end{equation}

In what follows  we will need the following technical result about the $L^1$-norm of the derivatives of the 
kernel $p(t,x,y)$,
\begin{lemma}\label{L1estimate}
If $p(t,x,y)$ is the Poisson-Hermite kernel,  then
\begin{equation}\label{1derivp-L^1-estimate}
\int_{\mathbb{R}^d} |\frac{\partial p(t,x,y)}{\partial t} | dy \leq \frac{C}{t}\, ,
\end{equation}
where $C$ is a constant independent of $x$ and $t$. Moreover, for any positive integer $k$ we have
\begin{equation} \label{kderivp-L^1-estimate}
\int_{\mathbb{R}^d} |\frac{\partial^k p(t,x,y)}{\partial t^k} | dy \leq \frac{C}{t^k}\, ,
\end{equation}
\end{lemma}
\dem
Let us  first  prove (\ref{1derivp-L^1-estimate}). Remember that
$$p(t,x,y) = \frac 1{2\pi ^{(d+1)/2}}\int_0^1t\frac{\exp \left(
t^2/4\log r\right) }{(-\log r)^{3/2}} \frac{\exp \left(
\frac{-\left| y-rx\right| ^2}{1-r^2}\right)
}{(1-r^2)^{d/2}}\frac{dr}r.$$
Therefore,
$$\frac{\partial p(t,x,y)}{\partial t} = \frac 1{2\pi ^{(d+1)/2}}\int_0^1\frac{\exp \left(
t^2/4\log r\right) }{(-\log r)^{3/2}}(1+ \frac{t^2}{2\log r}) \frac{\exp \left(
\frac{-\left| y-rx\right| ^2}{1-r^2}\right)
}{(1-r^2)^{d/2}}\frac{dr}r.$$
Then, by Tonelli's theorem, using  the fact that  
$$\frac{1}{\pi ^{d/2}}  \int_{\mathbb{R}^d} \frac{\exp \left(\frac{-\left| y-rx\right| ^2}{1-r^2}\right)}{(1-r^2)^{d/2}} dy =1,$$
\begin{eqnarray*}
\int_{\mathbb{R}^d} |\frac{\partial p(t,x,y)}{\partial t} | dy &\leq& \frac 1{2\pi ^{(d+1)/2}}\int_{\mathbb{R}^d}\int_0^1\frac{\exp \left(
t^2/4\log r\right) }{(-\log r)^{3/2}}|1+ \frac{t^2}{2\log r}| \frac{\exp \left(
\frac{-\left| y-rx\right| ^2}{1-r^2}\right)
}{(1-r^2)^{d/2}}\frac{dr}r \, dy \\
&=& \frac 1{2\pi ^{(d+1)/2}}\int_0^1 \frac{\exp \left(
t^2/4\log r\right) }{(-\log r)^{3/2}}|1+ \frac{t^2}{2\log r}|\int_{\mathbb{R}^d}  \frac{\exp \left(
\frac{-\left| y-rx\right| ^2}{1-r^2}\right)
}{(1-r^2)^{d/2}} \, dy \, \frac{dr}r\\
&=& \frac 1{2\pi ^{1/2}}\int_0^1 \frac{\exp \left(
t^2/4\log r\right) }{(-\log r)^{3/2}}|1+ \frac{t^2}{2\log r}| \, \frac{dr}r . 
\end{eqnarray*}
Thus what we need to prove is
\begin{equation}
\int_0^1 \frac{\exp \left(
t^2/4\log r\right) }{(-\log r)^{3/2}}|1+ \frac{t^2}{2\log r}| \, \frac{dr}r  \leq \frac{C}{t}.
\end{equation}
Making the change of variable $s= - \log r$ we get
\begin{eqnarray*}
\int_0^1 \frac{\exp \left(
t^2/4\log r\right) }{(-\log r)^{3/2}}|1+ \frac{t^2}{2\log r}| \, \frac{dr}r  &=& \int_0^{\infty} \frac{e^{
-t^2/4s} }{s^{3/2}}|1- \frac{t^2}{2s}| \, ds  \\
&\leq&\int_0^{\infty} \frac{e^{-t^2/4s} }{s^{3/2}}\, ds  + \int_0^{\infty} \frac{e^{-t^2/4s } }{s^{3/2}} \frac{t^2}{2s} \, ds  \\
\end{eqnarray*}
Now, making the change of variable $v= \frac{t^2}{4s}$, $ds =-\frac{t^2}{4v^2} dv$, we get
\begin{eqnarray*}
\int_0^{\infty} \frac{e^{-t^2/4s} }{s^{3/2}}\, ds  &=& \int_0^{\infty} e^{-v} (\frac{t^2}{4v})^{-3/2}\, \frac{t^2}{4v^2} dv\\
&=&  \int_0^{\infty} e^{-v} \frac{(4v)^{3/2}}{t^3}\, \frac{t^2}{4v^2} dv =\frac{C}{t} \int_0^{\infty} e^{-v} v^{-1/2} dv \\
&=&  \frac{C \Gamma(1/2)}{t} = \frac{C'}{t} 
\end{eqnarray*}
and 
\begin{eqnarray*}
\int_0^{\infty} \frac{e^{-t^2/4s} }{s^{3/2}} \frac{t^2}{4s}\, ds  &=& 2\int_0^{\infty} e^{-v} (\frac{t^2}{4v})^{-3/2} v\, \frac{t^2}{4v^2} dv\\
&=&  2\int_0^{\infty} e^{-v} \frac{(4v)^{3/2}}{t^3} v\, \frac{t^2}{4v^2} dv =\frac{C}{t} \int_0^{\infty} e^{-v} v^{1/2} dv \\
&=&  \frac{C \Gamma(3/2)}{t} = \frac{C'}{t}.
\end{eqnarray*}
For the proof of the general case (\ref{kderivp-L^1-estimate}) we use induction. Since the case $k=1$ is already proved let us assume that (\ref{kderivp-L^1-estimate}) holds for certain $k$ and prove that it also holds for $k+1$. By the semigroup property and taking $u=t+s$,  we have
\begin{eqnarray*}
\frac{\partial^{k+1} p(u,x,y)}{\partial u^{k+1}} &=&\frac{\partial}{\partial s} \frac{\partial^{k}}{\partial t^{k}} p(t+s,x,y)\\
&=& \frac{\partial}{\partial s} \frac{\partial^{k}}{\partial t^{k}} [ \int_{\mathbb{R}^d} p(s,x,v) p(t,v,y) dv] \\
&=&  \int_{\mathbb{R}^d} \frac{\partial p(s,x,v)}{\partial s} \frac{\partial^{k} p(t,v,y)}{\partial t^{k}} dv.\\
\end{eqnarray*}
Therefore, 
\begin{eqnarray*}
\int_{\mathbb{R}^d} |\frac{\partial^{k+1} p(u,x,y)}{\partial u^{k+1}} | dy &\leq& \int_{\mathbb{R}^d}  \int_{\mathbb{R}^d} |\frac{\partial p(s,x,v)}{\partial s}| |\frac{\partial^{k} p(t,v,y)}{\partial t^{k}}| dv dy\\
&\leq& \int_{\mathbb{R}^d} |\frac{\partial p(s,x,v)}{\partial s}| \int_{\mathbb{R}^d}|  \frac{\partial^{k} p(t,v,y)}{\partial t^{k}}| dy dv \\
&\leq&  \frac{C}{s}\frac{C}{t^k}.
\end{eqnarray*}
Finally,  taking $s=t= u/2$ the case $k+1$ is proved. \ep \\

As usual in what follows $C$ represents a constant that is not necessarily the same in each occurrence.\\

\section{Gaussian Lipschitz spaces.}
We want to define Lipschitz spaces associated to the Gaussian measure. First,  observe that   the spaces $L^p(\gamma)$ are not in general closed under the action of the classical  translation operator $\tau_yf(x) = f(x+y)$, see \cite{PinUr}. In the classical case, the Poisson semigroup provides an alternative characterization of the Lipchitz spaces, see \cite{se70}, we will follow this approach, using the Poisson-Hermite semigroup,  to define Gaussian Lipchitz spaces. First 
we need the following  key result,
 \begin{proposition}\label{Prop1} 
 Suppose $f \in L^{\infty}(\gamma)$ and $\alpha >0$. Let $k$ and $l$ be two integers both greater than $\alpha$. The two conditions
 \begin{equation}
\| \frac{\partial^k P_t f}{\partial t^k}\|_{\infty}  \leq A_{\alpha,k} t^{-k+\alpha}
\end{equation}
and 
 \begin{equation}
\| \frac{\partial^l P_t f}{\partial t^l}\|_{\infty}  \leq A_{\alpha,l} t^{-l+\alpha},
\end{equation}
are equivalent. Moreover, the smallest $A_{\alpha,k}$ and $A_{\alpha,l}$ holding in the above inequalities are comparable.
 \end{proposition}
\dem
It suffices to prove that if $k > \alpha$, 
 \begin{equation}\label{caso1}
\| \frac{\partial^k P_t f}{\partial t^k}\|_{\infty}  \leq A_{\alpha,k} t^{-k+\alpha}
\end{equation}
and 
 \begin{equation}\label{caso2}
\| \frac{\partial^{k+1} P_t f}{\partial t^{k+1}}\|_{\infty}  \leq A_{\alpha,k+1} t^{-(k+1)+\alpha},
\end{equation}
are equivalent.  \\

Let us assume (\ref{caso1}). Since by the semigroup property if $t=t_1+t_2$, $P_t f = P_{t_1} (P_{t_{2}} f)$ then  using the hypothesis and Lemma \ref{L1estimate},

\begin{eqnarray*}
 \| \frac{\partial^{k+1} P_t f}{\partial t^{k+1}} \|_{\infty} &=& \|  \frac{\partial P_{t_1}}{\partial t_1}(\frac{\partial^k P_{t_2} f}{\partial t_2^k})\|_{\infty}  \leq \| \frac{\partial^k P_{t_2} f}{\partial t_2^k}\|_{\infty} \int_{\mathbb{R}^d}| \frac{\partial p(t_1, \cdot, y)}{\partial t_1}| dy \\ 
 &\leq& A_{\alpha,k} t_2^{-k+\alpha} C t_1^{-1}.
\end{eqnarray*}
For $t_1 = t_2 = t/2$ we get (\ref{caso2}).\\

Let us now assume (\ref{caso2}). Observe that, again by Lemma \ref{L1estimate},
$$ \| \frac{\partial^k P_t f}{\partial t^k} \|_{\infty}  \leq  \|f\|_{\infty}  \int_{\mathbb{R}^d} |\frac{\partial^k p(t,x,y)}{\partial t^k} | dy  \leq \frac{C}{t^k} \|f\|_{\infty},$$
so $\frac{\partial^k P_t f}{\partial t^k} \rightarrow 0$ as $t \rightarrow \infty$, and then by hypothesis
$$  \| \frac{\partial^k P_t f}{\partial t^k} \|_{\infty}  \leq \int_t^{\infty} \| \frac{\partial^{k+1} P_s f}{\partial s^{k+1}} \|_{\infty}  ds \leq A_{\alpha,k+1} \frac{t^{-k+\alpha}}{-k+\alpha} = C t^{-k+\alpha}. $$ \ep \\

Now, we can define the Gaussian Lipschitz spaces as follows, 

\begin{defi} For  $ \alpha > 0$ let  $n$ be the smallest integer greater than $\alpha$. The Gaussian Lipschitz space $Lip_{\alpha}(\gamma)$ is defined as the set of  $L^{\infty}$ functions
 for which there exists a constant $A_{\alpha}(f)$ such that
\begin{equation}\label{poisson2}
\|\frac{\partial^{n}P_t f}{\partial
t^{n}}\|_{\infty}\leq A_{\alpha}(f) \,t^{-n+\alpha}.
\end{equation}
The norm of $f \in Lip_{\alpha}(\gamma)$ is defined as
\begin{equation}
\left\| f \right\|_{Lip_{\alpha}(\gamma)}: =  \left\| f
\right\|_{\infty} +A_{\alpha}(f),
\end{equation}
where $A_{\alpha}(f)$ is the smallest constant $A$ appearing in (\ref{poisson2}).
\end{defi}

\begin{obs} By Proposition \ref{Prop1} the definition of $Lip_{\alpha}(\gamma)$ does not depend on which  $k>\alpha$ is chosen and the resulting norms are equivalent. 
\end{obs}

\begin{obs} \label{steinobs} Condition  (\ref{poisson2}) is of interest for $t$ near zero, since the inequality
 \begin{equation}\label{poisson3}
\|\frac{\partial^{n}P_t f}{\partial
t^{n}}\|_{\infty}\leq At^{-n},
\end{equation}
which is stronger away from zero, follows for $f \in L^{\infty}$ immediately from (\ref{kderivp-L^1-estimate}),
\begin{eqnarray*}
\|\frac{\partial^{n}P_t f}{\partial
t^{n}}\|_{\infty} &\leq& \int_{\mathbb{R}^d} |\frac{\partial^n p(t,x,y)}{\partial t^n} | |f(y)| dy \leq \frac{C}{t^n} \| f\|_{\infty}. 
\end{eqnarray*}
\end{obs}

This spaces can be also obtained using abstract interpolation theory using the Poisson-Hermite semigroup,
see \cite{trie} 1.6.5.\\

Also there are inclusion relations among the Gaussian Lipschitz spaces,
 \begin{proposition}\label{Prop2} 
 If   $ 0 <\alpha_1 < \alpha_2 $ the inclusion  $Lip_{\alpha_2}(\gamma) \subset Lip_{\alpha_1}(\gamma)$ holds.
\end{proposition}

\dem
Take $f \in Lip_{\alpha_2}(\gamma)$ and $n \geq \alpha_2$, then 
$$ \|\frac{\partial^{n}P_t f}{\partial t^{n}}\|_{\infty}\leq A_{\alpha}(f)  t^{-n+\alpha_2}.$$
If $0 < t < 1$ then $t ^{-n+\alpha_2} \leq t^{-n+\alpha_1}$ and therefore
$$ \|\frac{\partial^{n}P_t f}{\partial
t^{n}}\|_{\infty}\leq A_{\alpha}(f)  t^{-n+\alpha_1}.$$
Now if $ t \geq 1$ then we know from (\ref{poisson3})  that 
$$ \|\frac{\partial^{n}P_t f}{\partial
t^{n}}\|_{\infty}\leq A_{\alpha}(f)  t^{-n} $$
and  as  $t ^{-n+\alpha_1} > t^{-n}$  we get also in this case
$$ \|\frac{\partial^{n}P_t f}{\partial
t^{n}}\|_{\infty}\leq A_{\alpha}(f) t^{-n+\alpha_1} $$
since $n > \alpha_1$; then $f \in Lip_{\alpha_1}(\gamma)$. 
\ep

\begin{proposition}\label{Prop3} 
If  $f \in  Lip_{\alpha}(\gamma)$, $0 <\alpha < 1$, then 
\begin{equation}\label{Lipcond}
|| P_t f - f||_{\infty} \leq  A_{\alpha}(f) \, t^{\alpha}.
\end{equation}
\end{proposition}
\dem
Applying the Fundamental Theorem of Calculus,
\begin{eqnarray*}
|| P_t f - f||_{\infty} &=& ||  \int_0^t \frac{\partial P_s f}{\partial
s} ds ||_{\infty} \leq \int_0^t  \|\frac{\partial P_s f}{\partial
s}\|_{\infty} ds\\
&\leq& A_{\alpha}(f)   \int_0^t  s^{-1+\alpha}ds = A_{\alpha}(f)  \, t^{\alpha}.
\end{eqnarray*}
\ep

\section{Boundedness of Fractional Integrals and Fractional Derivatives on $ Lip_{\alpha}(\gamma)$ }

 For $\beta>0$, the Bessel Potential of order $\beta>0,$
$\mathcal{J}_\beta$, associated to the Gaussian  measure
is defined formally as
\begin{eqnarray}
\mathcal{J}_\beta = (I+\sqrt{-L})^{-\beta},
\end{eqnarray}
meaning that for the Hermite polynomials we have,
\begin{eqnarray*}
\mathcal{J}_\beta h_\nu (x)=\frac 1{(1+\sqrt{\left|
\nu \right|})^{\beta}}h_\nu (x)
\end{eqnarray*}
By linearity this definition can be extended to any polynomial and P. A.  Meyer's theorem allows us to extend Bessel Potentials
to a \mbox{continuous} operator on $L^p(\gamma),$ $1 < p <
\infty$. Additionally, it is easy to see that $\mathcal{J}_\beta$ is a bijection over the set of polynomials ${\mathcal P}$. Alternatively the Bessel potentials can be
defined as
\begin{equation}\label{Beselrepre}
\mathcal{J}_\beta
f(x)=\frac{1}{\Gamma(\beta)}\int_{0}^{+\infty}s^{\beta-1}e^{-s}P_{s}f(x) ds,\; f \in L^p(\gamma).
\end{equation}
 For more details see \cite{forscotur}.
Moreover $\{\mathcal{J}_\beta\}_\beta$
is a strongly continuous semigroup on $L^p(\gamma)$, $1 \leq p <\infty$, with infinitesimal generator $\frac{1}{2}\log (I-L).$ \\

We will study the action of the Bessel potentials on the Gaussian Lipschitz spaces  $Lip_{\alpha}(\gamma)$

\begin{theorem}\label{continuidadBessel}
Let $\alpha> 0$ and $\beta>0$ then 
 $\mathcal{J}_{\beta}$ is
bounded from $Lip_{\alpha}(\gamma)$ to
 $Lip_{\alpha+\beta}(\gamma)$.
\end{theorem}
\dem

Let  $f\in Lip_{\alpha}(\gamma)$ and consider a fixed integer $n>\alpha+\beta$, then 
\begin{eqnarray*}
\|\frac{\partial^{n} P_{t}f}{\partial
t^{n}} \|_{\infty}&\leq&A_{\alpha}(f) t^{-n+\alpha},  \quad t>0.
\end{eqnarray*}

Using (\ref{Beselrepre}), the fact that $f \in L^{\infty}$, and consequently $P_{t+s} f \in L^{\infty}$,
 we obtain 

\begin{equation}\label{Beselrepre2}
P_{t}(\mathcal{J}^{\beta}f)(x)=\displaystyle\frac{1}{\Gamma(\beta)}\int_{0}^{+\infty}s^{\beta-1}e^{-s}P_{t+s}f(x)\, ds,
\end{equation}
and therefore 
$$ \| P_{t}(\mathcal{J}_{\beta}f)\|_{\infty} \leq \| f\|_{\infty},$$
i.e. $P_{t}(\mathcal{J}_{\beta}f) \in L^{\infty}$. \\

Now we want to verify the Lipchitz condition. Differentiating  (\ref{Beselrepre2}), we get  

\begin{eqnarray*}
\frac{\partial^{n}P_{t}(\mathcal{J}_{\beta}f)(x)}{\partial t^{n}}&=&\frac{1}{\Gamma(\beta)}\int_{0}^{+\infty}s^{\beta-1}e^{-s}\frac{\partial^{n} P_{t+s}f(x)}{\partial t^{n}}\, ds\\
&=&\frac{1}{\Gamma(\beta)}\int_{0}^{+\infty}s^{\beta-1}e^{-s}\frac{\partial^{n}
P_{t+s}f(x)}{\partial (t+s)^{n}}\, ds,
\end{eqnarray*}
and this implies

\begin{eqnarray*}
\|\frac{\partial^{n}P_{t}(\mathcal{J}_{\beta}f)}{\partial t^{n}}\|_{\infty}&\leq& \frac{1}{\Gamma(\beta)}\int_{0}^{t}s^{\beta-1}e^{-s}\|\frac{\partial^{n}P_{t+s}f}{\partial
(t+s)^{n}}\|_{\infty}\, ds \\
& & \quad \quad \quad +  \frac{1}{\Gamma(\beta)}\int_{t}^{+\infty}s^{\beta-1}e^{-s}\|\frac{\partial^{n}P_{t+s}f}{\partial (t+s)^{n}}\|_{\infty}\, ds\\
&=& (I)+(II).
\end{eqnarray*}
Since $\beta>0$ and $t+s>t $,
\begin{eqnarray*}
(I)&\leq&\frac{A_{\alpha}(f) }{\Gamma(\beta)}\int_{0}^{t}s^{\beta-1} (t+s)^{-n+\alpha} e^{-s}\, ds\\
&\leq&\frac{A_{\alpha}(f) }{\Gamma(\beta)} t^{-n+\alpha} \int_{0}^{t}s^{\beta-1}ds (\gamma) \leq  C t^{-n+\alpha+\beta} \, \| f \|_{Lip_{\alpha}(\gamma)}.
\end{eqnarray*}

On the other hand, since $n>\alpha+\beta$, and $t+s>s$ 
\begin{eqnarray*}
(II) &\leq&\frac{A_{\alpha}(f) }{\Gamma(\beta)}\int_{t}^{\infty}s^{\beta-1}e^{-s}(t+s)^{-n+\alpha}\, ds\\
&\leq&\frac{A_{\alpha}(f) }{\Gamma(\beta)}\int_{t}^{\infty}s^{\beta-1}e^{-s} s^{-n+\alpha}\, ds\\
&\leq&\frac{ A_{\alpha}(f) }{\Gamma(\beta)}\int_{t}^{\infty}s^{-n+\alpha+\beta-1}ds=C  A_{\alpha}(f) t^{-n+\alpha+\beta}.
\end{eqnarray*}
Therefore,
\begin{eqnarray*}
\|\frac{\partial^{n}P_{t}(\mathcal{J}_{\beta}f)}{\partial t^{n}}\|_{\infty}&\leq&C  A_{\alpha}(f) t^{-n+\alpha+\beta},  \quad t>0.
\end{eqnarray*} 
Thus
$\mathcal{J}_{\beta}f\in Lip_{\alpha+\beta}(\gamma)$ and moreover, 
\begin{eqnarray*}
\|{\mathcal{J}}_{\beta}f\|_{Lip_{\alpha+\beta}(\gamma)}&=&\|{\mathcal{J}}_{\beta}f\|_{\infty}+A_{\alpha}({\mathcal{J}}_{\beta}f)\\
&\leq&\|f\|_{\infty}+CA_{\alpha}(f)  \leq C \|f\|_{Lip_{\alpha}(\gamma)}. 
\end{eqnarray*} 
 \ep
 
For $\beta>0$ the Riesz fractional integral or Riesz potential of order $\beta$, $I_\beta$, with respect to the Gaussian measure is defined formally as
\begin{equation}\label{i1}
I_\beta=(-L)^{-\beta/2}\Pi_{0},
\end{equation}
where,
$\Pi_{0}f=f-\displaystyle\int_{\mathbb{R}^{d}}f(y)\gamma(dy)$,
for $f\in L^{2}(\gamma)$. That means that for  the Hermite
polynomials $\{h_\nu\}$, with $\left|\nu\right|>0$,
\begin{equation}\label{e4}
I_\beta h_\nu(x)=\frac 1{\left|
\nu\right|^{\beta/2}}h_\nu(x),
\end{equation}
and for $\nu=(0,0, \cdots,0), \,
I_\beta (h_{\nu})=0.$ By linearity can be extended to any polynomial.
If $f$ is a polynomial with $ \int_{\mathbb{R}^d} f(y) \gamma(dy) = 0,$
\begin{equation}\label{e3}
I_\beta f(x)  =\frac 1{\Gamma(\beta)}\int_0^{\infty}
s^{\beta-1}P_s f(x)\,ds.
\end{equation}
By P. A. Meyer's multiplier theorem, $I_\beta$ admits a continuous extension  to $ L^p(\gamma _d)$,  $1 < p < \infty$, and (\ref{e3}) can be extended for $f \in L^{p}(\gamma)$, see \cite{sp97}. In addition if $f$ $\in C_B^2(\mathbb{R}^d)$ such that $\int_{\mathbb{R}^d}f(y)\gamma(dy)=0$, then
\begin{equation}\label{p2}
I_\beta f =-\frac {1}{\beta \Gamma(\beta)}\int_0^{\infty} s^\beta \frac{\partial P_s f }{\partial s} ds,
\end{equation}
see \cite{lour}.

The Riesz Potentials are not bounded operators on $ L^{\infty}(\gamma _d)$ nor on  $Lip_{\alpha} (\gamma)$.\\ 

Following the classical case, the Riesz fractional derivate of order $\beta >0$ with respect to the Gaussian 
measure $D^\beta$, is defined formally as
\[
D^\beta=(-L)^{\beta/2},
\]
which means that for the Hermite polynomials, we have
\begin{equation}\label{e6}
D^\beta h_\nu(x)=\left| \nu\right|^{\beta/2} h_\nu(x),
\end{equation}
In the case of $0 < \beta < \alpha <1$ we have the following integral representation,  
\begin{equation}\label{e5}
D^\beta f =\frac 1{c_\beta}\int_0^{\infty}s^{-\beta-1}(P_s  -I)\, f ds,
\end{equation} 
for $f \in Lip_{\alpha} (\gamma)$, where
$
c_\beta=\int_0^\infty u^{-\beta-1}(e^{-u}-1)du,$  which is finite as  $0 < \beta <1.$\\

The action of  Riesz fractional derivates on Lipchitz spaces is as follows,

\begin{theorem}\label{actfracdev} For $0<\beta < \alpha <1$,  the Riesz fractional derivate of order $\beta$, $ D^\beta: Lip_{\alpha}(\gamma)  \rightarrow Lip_{\alpha-\beta}(\gamma)$ is bounded.
\end{theorem}
\dem
Take $f \in Lip_{\alpha}(\gamma)$ i.e. $f\in L^{\infty}$ and $\|\frac{\partial P_t f}{\partial
t}\|_{\infty}\leq A_{\alpha}(f) t^{-1+\alpha}$.   Let us observe that using representation (\ref{e5}), Proposition \ref{Prop1} and (\ref{Lipcond}) we get,
\begin{eqnarray*}
| D^\beta f(x)| &\leq& \frac {1}{ c_{\beta}}\int_0^1 s^{-\beta-1}| P_s f(x) -f(x) | ds+ \frac {1}{ c_{\beta}}\int_1^{\infty} s^{-\beta-1}|P_s f(x) -f(x)| ds\\
&\leq&\frac {1}{c_{\beta}}\int_0^1 s^{-\beta-1} \| P_s f -f \|_{\infty}  ds+ \frac {2 \|f\|_{\infty}  }{ c_{\beta}}\int_1^{\infty} s^{-\beta-1}ds\\
&\leq&\frac {A_{\alpha}(f) }{ c_{\beta}}\int_0^1 s^{\alpha -\beta-1} ds+ \frac {2 \|f\|_{\infty, \gamma}   }{ c_{\beta}}\int_1^{\infty} s^{-\beta-1}ds \\
&=& \frac {A_{\alpha}(f) }{ c_{\beta}(\alpha -\beta)}+ \frac {2 \|f\|_{\infty}   }{\beta c_{\beta}} \leq C_{\alpha, \beta} \| f \|_{Lip_{\alpha}(\gamma)},
\end{eqnarray*}
thus $D^\beta f \in L^{\infty}$. \\

Now we want to verify the Lipchitz condition. Fixing $t$ and using again representation (\ref{e5}), we have
\begin{eqnarray*}
\frac{\partial}{\partial t}(P_t  D^{\beta} f(x)) &=&  \frac {1}{ c_{\beta}}\frac{\partial}{\partial t}[\int_0^{\infty} s^{-\beta-1} (P_{t+s} f(x) -P_t f(x)) ds]\\
&=&  \frac {1}{ c_{\beta}}\int_0^{\infty} s^{-\beta-1} [\frac{\partial P_{s+t} f(x)}{\partial t} -\frac{\partial P_t f(x)}{\partial s}  ] ds\\
&=&  \frac {1}{ c_{\beta}}\int_0^{t} s^{-\beta-1} [\frac{\partial P_{s+t} f(x)}{\partial t} -\frac{\partial P_t f(x)}{\partial t}  ] ds \\
&&  \quad +\frac {1}{ c_{\beta}}\int_t^{\infty} s^{-\beta-1} [\frac{\partial P_{s+t} f(x)}{\partial t} -\frac{\partial P_t f(x)}{\partial t}  ] ds\\
&=& (I) + (II).
\end{eqnarray*}
By Proposition \ref{Prop2}  we have
\begin{equation}\label{estimate2}
\|\frac{\partial^2 P_r f}{\partial r^2}  \|_{\infty}  \leq A\, r^{\alpha -2},
\end{equation}
and using the Fundamental Theorem of Calculus we get, for $s >0$
\begin{eqnarray}\label{TFC2}
\nonumber | \frac{\partial P_{t+s} f(x)}{\partial t} -\frac{\partial  P_t f(x)}{\partial t}| &\leq& \int_t^{t+s} |\frac{\partial^2  P_r f(x)}{\partial r^2}| dr 
\leq A \int_t^{s+t}  r^{\alpha  -2} dr\\
& \leq& \, A t^{\alpha  -2} s.
\end{eqnarray}

Therefore, 
\begin{eqnarray*}
 |(I)| &\leq& \frac {1}{ c_{\beta}}\int_0^{t} s^{-\beta-1} |\frac{\partial P_{t+s} f(x)}{\partial t} -\frac{\partial  P_t f(x)}{\partial t} | ds\\
 &\leq& A \frac {t^{\alpha-2}}{ c_{\beta}}\int_0^{t} s^{-\beta} ds = C_{\alpha, \beta}  \, t^{-1+\alpha-\beta}.
\end{eqnarray*}
On the other hand, using (\ref{poisson3}),
\begin{eqnarray*}
 |(II)| &\leq& \frac {1}{ c_{\beta}}\int_t^{\infty} s^{-\beta-1} [|\frac{\partial P_{t+s} f(x)}{\partial t}| +|\frac{\partial  P_t f(x)}{\partial t}  |] ds\\
 &\leq& \frac {C}{ c_{\beta}}\int_t^{\infty} s^{-\beta-1} [(t+s)^{-1+\alpha}+t^{-1+\alpha}] ds\\
 &\leq& Ct^{-1+\alpha} \int_t^{\infty} s^{-\beta-1} ds = C_{\alpha, \beta}  \,t^{-1+\alpha-\beta}.
\end{eqnarray*}

Thus,
$$ \| \frac{\partial}{\partial t}(P_t  D^\beta f)\|_{\infty}  \leq C_{\alpha, \beta} t^{\alpha-\beta-1}, $$
which implies $D^\beta f\in Lip_{\alpha-\beta}(\gamma).$
\ep \\

We can also consider a Bessel fractional derivative $ {\mathcal D}^\beta$,  defined formally as
\[
 {\mathcal D}^\beta=(I+\sqrt{-L})^{\beta},
\]
which means that for the Hermite polynomials, we have
\begin{equation}\label{e7}
 {\mathcal D}^\beta h_\nu(x)=(1+ \sqrt{\left| \nu\right|})^{\beta} h_\nu(x),
\end{equation}
In the case of $0 < \beta <1$ we have the following integral representation,  
\begin{equation}\label{e8}
 {\mathcal D}^\beta f =\frac 1{c_\beta}\int_0^{\infty}t^{-\beta-1}( e^{-t} P_t -I) \, f dt,
\end{equation} 
where, as before, 
$
c_\beta=\int_0^\infty u^{-\beta-1}(e^{-u}-1)du.$\\

We want to study of the  action of the Bessel fractional derivative ${\mathcal D}^\beta$ on the Gaussian Lipschitz spaces.

\begin{theorem} For $0<\beta < \alpha <1$,  the Bessel fractional derivate of order $\beta$, ${\mathcal D}^\beta: Lip_{\alpha}(\gamma)  \rightarrow Lip_{\alpha-\beta}(\gamma)$ is bounded.
\end{theorem}
\dem

The proof of this result is essentially analogous to the proof of Theorem \ref{actfracdev}.
Let  $f \in Lip_{\alpha}(\gamma)$ i.e. $f\in L^{\infty}$ such that $\|\frac{\partial P_t f}{\partial
t}\|_{\infty}\leq A_{\alpha}(f) t^{-1+\alpha}$.   Using the representation (\ref{e8}), (\ref{Lipcond}) and Proposition \ref{Prop1}, we get,
\begin{eqnarray*}
| {\mathcal D}^\beta f(x)| &\leq& \frac {1}{ c_{\beta}}\int_0^{\infty} s^{-\beta-1}| e^{-s}P_s f(x) -f(x) | ds\\
&=&\frac {1}{ c_{\beta}}\int_0^1 s^{-\beta-1}| e^{-s} P_s f(x) -f(x) | ds+ \frac {1}{ c_{\beta}}\int_1^{\infty} s^{-\beta-1}|e^{-s}P_s f(x) -f(x)| ds\\
&=& (I) + (II).
\end{eqnarray*}
Now 
\begin{eqnarray*}
(I) &\leq& \frac {1}{c_{\beta}}\int_0^1 s^{-\beta-1} e^{-s} |  P_s f(x)  - f(x)  |ds +  \frac {1}{c_{\beta}}\int_0^1 s^{-\beta-1} |e^{-s} -1 | \, |f(x)| ds \\
&\leq& \frac {1}{c_{\beta}}\int_0^1 s^{-\beta-1}  \| P_s f -f \|_{\infty}\, ds +  \frac {1}{c_{\beta}}\int_0^1 s^{-\beta-1} |e^{-s} -1 | \, |f(x)| ds \\
&\leq&\frac {A_{\alpha}(f) }{ c_{\beta}}\int_0^1 s^{\alpha -\beta-1} ds+ \frac { \|f\|_{\infty}   }{ c_{\beta}}\int_0^1 s^{-\beta-1} (1- e^{-s}) ds \\
&=& \frac {A_{\alpha}(f) }{ c_{\beta}(\alpha -\beta)}+ \frac { C \|f\|_{\infty}   }{ c_{\beta}} \leq C_{\alpha, \beta} \| f \|_{Lip_{\alpha}(\gamma)}.
\end{eqnarray*}
and \begin{eqnarray*}
(II) &\leq&  \frac {1}{ c_{\beta}}\int_1^{\infty} s^{-\beta-1}[e^{-s}|P_s f(x)| +|f(x)| ]ds\\
&\leq & \frac {2 \|f\|_{\infty}  }{ c_{\beta}}\int_1^{\infty} s^{-\beta-1}ds \leq \frac {2 \|f\|_{\infty}   }{\beta c_{\beta}}
\leq C_{\alpha, \beta} \| f \|_{Lip_{\alpha}(\gamma)},
\end{eqnarray*}
thus $  {\mathcal D}^\beta f \in L^{\infty}$.\\

Now we want to verify the Lipchitz condition. By Observation \ref{steinobs} it is enough to consider the case $0<t<1/2$. Using again representation (\ref{e8}) we have
\begin{eqnarray*}
\frac{\partial}{\partial t}(P_t   {\mathcal D}^\beta f(x)) &=&  \frac {1}{ c_{\beta}}\frac{\partial}{\partial t}[\int_0^{\infty} s^{-\beta-1} ( e^{-s} P_{t+s} f(x) -P_t f(x) ) ds]\\
&=&  \frac {1}{ c_{\beta}}\int_0^{\infty} s^{-\beta-1} [ e^{-s} \frac{\partial P_{t+s} f(x)}{\partial t} - \frac{\partial P_{t} f(x)}{\partial t} ] ds\\
&=&   \frac {1}{ c_{\beta}}\int_0^{t} s^{-\beta-1} [ e^{-s} \frac{\partial P_{t+s} f(x)}{\partial t} - \frac{\partial P_{t} f(x)}{\partial t} ] ds\\
&&  \quad + \frac {1}{ c_{\beta}}\int_t^{\infty} s^{-\beta-1} [ e^{-s} \frac{\partial P_{t+s} f(x)}{\partial t} - \frac{\partial P_{t} f(x)}{\partial t} ] ds\\
&=& (III)+(IV).
\end{eqnarray*}
 Using (\ref{estimate2}) we have, for $ 0< r< 1,$
\begin{eqnarray*}
|\frac{\partial [e^{-r}  \frac{\partial P_r f(x)}{\partial r}]}{\partial r}  |  &\leq& |-e^{-r}  \frac{\partial P_r f(x)}{\partial r}  |   +  | e^{-r}  \frac{\partial^2 P_r f(x)}{\partial r^2}  |  \\
&\leq&  A_{\alpha}(f)  e^{-r} \, r^{\alpha -1}+ A e^{-r} \, r^{\alpha -2} < C e^{-r} \, r^{\alpha -2},
\end{eqnarray*}
and then, by the Fundamental Theorem of Calculus, for $s< t$
\begin{eqnarray}\label{EstB1}
\nonumber  |e^{-(t+s)} \frac{\partial P_{t+s} f(x)}{\partial t} - e^{-t}\frac{\partial P_{t} f(x)}{\partial t}|  &\leq&
\int_t^{t+s}   |\frac{\partial [e^{-r}  \frac{\partial P_r f(x)}{\partial r}]}{\partial r}  | dr \\
&\leq& C  \int_t^{t+s} e^{-r} \, r^{\alpha -2} dr \leq C e^{-t} \,  t^{\alpha -2} \, s. 
\end{eqnarray}
Thus
 \begin{eqnarray*}
 |(III)|  &\leq&  \frac {e^t}{ c_{\beta}}\int_0^{t} s^{-\beta-1} |e^{-(t+s)} \frac{\partial P_{t+s} f(x)}{\partial t} - e^{-t}\frac{\partial P_{t} f(x)}{\partial t}| ds\\
  &\leq& \frac {C e^t}{ c_{\beta}}\int_0^{t} s^{-\beta-1} e^{-t} \, t^{\alpha -2} s \, ds = C t^{\alpha-2} \int_0^{t} s^{-\beta} ds = C_{\alpha, \beta}  \, t^{-1 +\alpha-\beta}.
\end{eqnarray*}

On the other hand,  as $t+s > t$ we have
\begin{eqnarray*}
|(IV)|&\leq& \frac {1}{ c_{\beta}}\int_t^{\infty} s^{-\beta-1} | e^{-s} \frac{\partial P_{t+s} f(x)}{\partial t} - \frac{\partial P_{t} f(x)}{\partial t} | ds\\
&\leq&\frac {1}{ c_{\beta}}\int_t^{\infty} s^{-\beta-1} [|\frac{\partial P_{t+s} f(x)}{\partial t}| +|\frac{\partial  P_t f(x)}{\partial t}  |] ds\\
 &\leq& \frac {C}{ c_{\beta}}\int_t^{\infty} s^{-\beta-1} [(t+s)^{-1+\alpha}+t^{-1+\alpha}] ds\\
 &\leq& Ct^{\alpha-1} \int_t^{\infty} s^{-\beta-1} ds = C_{\alpha, \beta}  \,t^{-1+\alpha-\beta}.
\end{eqnarray*}

Therefore, 
$$ \| \frac{\partial}{\partial t}(P_t  {\mathcal D}^\beta f)\|_{\infty}  \leq C_{\alpha, \beta} t^{-1+\alpha-\beta}, $$
which implies $ {\mathcal D}^\beta f\in Lip_{\alpha-\beta}(\gamma).$ \ep \\
 
 Moreover, if $\beta \geq 1$, let  $k$ be  the smallest integer such that $ \beta < k$, then the Riesz Fractional Derivative $D^\beta$ can be represented as 
\begin{equation}\label{kRder}
D^\beta f = \frac{1}{c^k_{\beta}}\int_0^{\infty} s^{-\beta-1} ( P_s -I )^k f \, ds,
\end{equation}
 and the Bessel fractional derivative ${\mathcal D}^\beta$ can be represented  as 
\begin{equation}\label{kBder}
{\mathcal D}^\beta f =  \frac{1}{c^k_{\beta}} \int_0^{\infty} s^{-\beta-1} (e^{-s} P_s-I )^k\,  f \, ds,
\end{equation}
where in both cases $c^k_{\beta} =  \int_0^{\infty} u^{-\beta-1} (e^{-u} -1 )^k du < \infty$.

Observe that  (\ref{kRder}) and (\ref{kBder}) are the right formulas since it is easy to prove form  that for any Hermite polynomial $h_{\nu}$, 
$$ D^\beta h_{\nu} = \nu^{\beta/2} h_{\nu}, \quad \mbox{and} \quad {\mathcal D}^\beta h_{\nu} =(1+\sqrt{ \nu})^{\beta} h_{\nu}.$$

In this general case we also want to study of the  action of $D^\beta$ and ${\mathcal D}^\beta$ on the Gaussian Lipschitz spaces,
\begin{theorem}\label{actfracdevk} Given  $1 \leq \beta < \alpha $,  then
\begin{enumerate}
\item [i)] The Riesz fractional derivate of order $\beta$, $ D^\beta: Lip_{\alpha}(\gamma)  \rightarrow Lip_{\alpha-\beta}(\gamma)$ is bounded.
\item [ii)] The Bessel fractional derivate of order $\beta$, ${\mathcal D}^\beta: Lip_{\alpha}(\gamma)  \rightarrow Lip_{\alpha-\beta}(\gamma)$ is bounded.
\end{enumerate}
\end{theorem}

 First of all, observe that, using the Binomial Theorem and the semigroup property, we have 
 \begin{eqnarray}\label{powerrep}
\nonumber ( P_t -I )^k f(x) &=& \sum_{j=0}^k {k \choose j} P_t^{k-j} (-I)^j f(x) = \sum_{j=0}^k {k \choose j} (-1)^jP_t^{k-j} f(x)\\
\nonumber &=&\sum_{j=0}^k {k \choose j} (-1)^jP_{(k-j)t} f(x) =\sum_{j=0}^k {k \choose j} (-1)^j u(x,(k-j)t)\\
&=& \Delta_t^k (u(x, \cdot), 0),
\end{eqnarray}
where  as usual, $u(x,t) = P_t f(x)$, and 
$$\Delta_s^k (f,t)= \sum_{j=0}^k {k \choose j} (-1)^j f( t+(k-j)s)$$
is the $k$-th order forward difference of $f$ starting at $t$ with increment $s$.
We will need some technical results about forward differences that will be used later. These are well known results in forward differences' theory, see for instance \cite{fort}, but for the sake of completeness, their proofs  will be given in an appendix. 

\begin{lemma} \label{lemaforwdiff}
The forward differences have the following properties,
\begin{enumerate}
\item[i)] For any positive integer $k$, 
\begin{equation}\label{iterat}
\Delta^k_s(f,t) = \Delta_s^{k-1}(\Delta_s(f,\cdot),t) = \Delta_s(\Delta_s^{k-1}(f,\cdot),t). 
\end{equation}
\item[ii)] For any positive integer $k$,
\begin{equation}\label{relDifeDer}
\Delta^k_s(f,t) = \int_t^{t+s} \int_{v_1}^{v_1+s}  \cdots \int_{v_{k-2}}^{v_{k-2}+s}\int_{v_{k-1}}^{v_{k-1}+s} f^{(k)}(v_k) \, dv_k dv_{k-1} \cdots dv_2 dv_1
\end{equation}
\item [iii)] For any positive integer $k$,
\begin{equation}\label{difder}
\frac{\partial }{\partial s}(\Delta_s^k (f,t) ) = k \,\Delta_s^{k-1} (f',t+s),
\end{equation}
and for any integer $j>0$,
\begin{equation}\label{difder2}
\frac{\partial^j }{\partial t^j}(\Delta_s^k (f,t) ) =\Delta_s^{k} (f^{(j)},t).
\end{equation}
\end{enumerate}
\end{lemma}

For the proof  of Theorem \ref{actfracdevk} we will need  estimates analogous to (\ref{Lipcond}), (\ref{TFC2}), and those will  follow from the next result.
\begin{proposition} Let $\delta$ a real number and $k $ a positive integer  such that $ \delta <k$. Let $f$ be a function such that for some integer $k$
\begin{equation}\label{condder}
|f^{(k)}(r) | \leq C r^{-k+\delta},
\end{equation}
then
\begin{equation}\label{difder1}
|\Delta_s^k (f,t)| \leq C\, s^k t^{-k+\delta}, 
\end{equation}
\end{proposition}

\dem
The proof is inmmediate from (\ref{relDifeDer}), since as $ \delta <k$ 
\begin{eqnarray*}
|\Delta^k_s(f,t)| &\leq& \int_t^{t+s} \int_{v_1}^{v_1+s}  \cdots \int_{v_{k-2}}^{v_{k-2}+s}\int_{v_{k-1}}^{v_{k-1}+s} |f^{(k)}(v_k)| \, dv_k dv_{k-1} \cdots dv_2 dv_1\\
&\leq&C \int_t^{t+s} \int_{v_1}^{v_1+s}  \cdots \int_{v_{k-2}}^{v_{k-2}+s}\int_{v_{k-1}}^{v_{k-1}+s} v_k^{-k+\delta} \, dv_k dv_{k-1} \cdots dv_2 dv_1\\
&\leq&C \int_t^{t+s} \int_{v_1}^{v_1+s} \cdots \int_{v_{k-2}}^{v_{k-2}+s} s v_{k-1}^{-k+\delta} \, dv_{k-1} \cdots dv_2 dv_1\\
& & \cdots \cdots  \cdots\\
&\leq& C s^k t^{-k+\delta}.
\end{eqnarray*}
\ep

The following result is a generalization of Proposition \ref{Prop3}, 
\begin{proposition} \label{Propk} We have the following estimates,
\begin{enumerate}
\item[i)] If $f \in L^{\infty}$, for any positive integer $k$
  \begin{equation} \label{Poissboundn}
|| (P_t  - I)^k f||_{\infty} \leq  2^k\,\|f\|_{\infty}
\end{equation}
\item[ii)] Let $\alpha > 1$ and $n$ be the smallest integer bigger than $\alpha$. If $f \in Lip_{\alpha} (\gamma)$ then 
  \begin{equation} \label{Lipcond2}
|| (P_t  - I)^n f||_{\infty} \leq  A_{\alpha}(f) \, t^{\alpha}.
\end{equation}
\end{enumerate}
\end{proposition}

\dem

\begin{enumerate}
\item[i)] We already know from (\ref{powerrep}) that
$$ ( P_t -I )^k f(x) =  \Delta_t^n (u(x, \cdot), 0).$$ 
Then for any $k$  inequality (\ref{Poissboundn}) is immediate,
\begin{eqnarray*}
|| (P_t  - I)^k f||_{\infty} \leq  \sum_{j=0}^k {k \choose j}\| P_{(k-j)t} f\|_{\infty} = 2^k \|f\|_{\infty}.
\end{eqnarray*}

\item[ii)] Now to prove (\ref{Lipcond2}) observe  $\alpha -1 < n-1$ and condition (\ref{poisson2}) can be rewritten as
$$ \| \frac{\partial^{n}}{\partial t^{n}} (u(\cdot, t))\|_{\infty} = \| \frac{\partial^{n-1}}{\partial t^{n-1}} (u'(\cdot, t))\|_{\infty} \leq A_{\alpha}(f)  t^{-n+1+(\alpha-1)},$$
i.e condition (\ref{condder}) is satisfied for $\delta= \alpha-1$, then using (\ref{difder}) and then (\ref{difder2}) with $t=s=r$, 
\begin{eqnarray*}
 |( P_t -I )^n f(x)| &\leq & \int_0^t | \frac{\partial}{\partial r} ( \Delta_r^n (u(x, \cdot), 0))| dr\\
 &=& n \int_0^t |( \Delta_r^{n-1} (u'(x, \cdot), r)| dr\\
  &\leq& n A_{\alpha}(f)  \int_0^t  r^{n-1} r^{-n+1+(\alpha-1)}dr = C \int_0^t  r^{\alpha-1} dr = C t^{\alpha}. 
\end{eqnarray*}
\end{enumerate}
\ep

Finally, let us prove Theorem \ref{actfracdevk}.\\
\begin{enumerate}
\item [i)] Take $f \in Lip_{\alpha}(\gamma)$, then $f\in L^{\infty}$ and $\|\frac{\partial^n u(\cdot,t)}{\partial
t^n}\|_{\infty}\leq A_{\alpha}(f) t^{-n+\alpha}$.  Remember that $\beta < \alpha$,  $k$ is  the smallest integer bigger than $ \beta$ and let $n$ be the smallest integer bigger than $\alpha$, note $k \leq n$. Using representation (\ref{kRder}), and then inequalities (\ref{Poissboundn}) and (\ref{Lipcond2}),
\begin{eqnarray*}
| D^\beta f(x)| &\leq& \frac{1}{c^k_{\beta}}\int_0^{\infty} s^{-\beta-1} |( P_s -I )^k f(x)| \, ds\\ 
&=&\frac {1}{ c^k_{\beta}}\int_0^1 s^{-\beta-1} |( P_s -I )^k f(x)| ds+ \frac {1}{ c^k_{\beta}}\int_1^{\infty} s^{-\beta-1} |( P_s -I )^k f(x)| ds\\
&=& (I)+ (II).
\end{eqnarray*}
Now, let us assume $k<n$, if $k=n$ the argument is straightforward. Let $\varepsilon >0$ such that $\beta+\varepsilon <k$, by Proposition \ref{Prop2} $Lip_{\alpha}(\gamma) \subset Lip_{\beta+\varepsilon}(\gamma)$,  using (\ref{Lipcond2}),
\begin{eqnarray*}
(I) &\leq& \frac {1}{c^k_{\beta}}\int_0^1 s^{-\beta-1}  \| (P_s - I)^k f \|_{\infty}\, ds \leq \frac {A_{\beta+\varepsilon}(f) }{ c^k_{\beta}}\int_0^1 s^{\beta+\varepsilon -\beta-1} ds = \frac {A_{\beta+\varepsilon}(f) }{ c^k_{\beta}\varepsilon}.
\end{eqnarray*}
On the other hand
\begin{eqnarray*}
(II) &\leq&  \frac {1}{c^k_{\beta}}\int_1^{\infty} s^{-\beta-1}  \| (P_s - I)^k f \|_{\infty} \,  ds \leq \frac {2^k \|f\|_{\infty}   }{ c^k_{\beta}}\int_1^{\infty} s^{-\beta-1}  ds =  C_{\beta} \|f\|_{\infty}.
\end{eqnarray*}
Thus  $ D^\beta f \in L^{\infty}$.\\

Now we want to verify the Lipchitz condition. Observe that by the semigroup property,
\begin{eqnarray*}
P_t [( P_s -I )^k f(x) ]  &=& P_t (\Delta_s^k (u(x, \cdot), 0)) = P_t ( \sum_{j=0}^k {k \choose j} (-1)^jP_{(k-j)s} f(x) )\\
&=&    \sum_{j=0}^k {k \choose j} (-1)^jP_{t+(k-j)s} f(x) = \Delta_s^k (u(x, \cdot), t).
\end{eqnarray*}
Fixing $t>0$, using  again representation (\ref{kRder}) and (\ref{difder2}), we have 
\begin{eqnarray*}
\frac{\partial^n(P_t  D^\beta f(x))}{\partial t^n} &=&  \frac {1}{ c^k_{\beta}}\frac{\partial^n}{\partial t^n} [\int_0^{\infty} s^{-\beta-1} P_t [( P_s -I )^k f(x) ] ds]\\
&=&  \frac {1}{ c^k_{\beta}}\int_0^{\infty} s^{-\beta-1}\frac{\partial^n}{\partial t^n}[ \Delta^k_s(u(x,\cdot),t)] ds\\
&=&  \frac {1}{ c^n_{\beta}}\int_0^{t} s^{-\beta-1}  [ \Delta^k_s(u^{(n)}(x,\cdot),t)]  ds \\
&&  \quad +\frac {1}{ c^n_{\beta}}\int_t^{\infty} s^{-\beta-1} [ \Delta^k_s(u^{(n)}(x,\cdot),t)]  ds\\
&=& (III) + (IV).
\end{eqnarray*}

Now, by Proposition  \ref{Prop1}  we have from (\ref{poisson2}),
$$ \| \frac{\partial^{k}}{\partial t^{k}} (u^{(n)}(\cdot, t))\|_{\infty} =  \| \frac{\partial^{n+k} (u(\cdot, t))}{\partial t^{n+k}} \|_{\infty}  \leq A t^{-(n+k) +\alpha}= A t^{-k +(\alpha-n)},$$
 then by (\ref{difder1})  
\begin{eqnarray*}
 |(III)| &\leq& \frac {1}{ c^k_{\beta}}\int_0^{t} s^{-\beta-1} |  \Delta^k_s(u^{(n)}(x,\cdot),t) | ds\\
 &\leq&  \frac { A\, t^{-k+(\alpha-n)}}{ c^k_{\beta}}\int_0^{t} s^{-\beta+k-1} ds = C_{\alpha, \beta}  \, t^{-k+\alpha-n} t^{-\beta+k} = C_{\alpha, \beta}  \, t^{-n+\alpha-\beta}.
\end{eqnarray*}
On the other hand, 
\begin{eqnarray*}
|  \Delta^k_s(u^{(n)}(x,\cdot),t) | &\leq& \sum_{j=0}^k {k \choose j} |u^{(n)}(x,t+(k-j)s)|\\
&\leq& A_{\alpha}(f) \sum_{j=0}^k {k \choose j} |(t+(k-j)s)^{-n+\alpha}| \leq C t^{-n+\alpha}
\end{eqnarray*}
and then,
\begin{eqnarray*}
 |(IV)| &\leq& \frac {1}{ c^n_{\beta}}\int_t^{\infty} s^{-\beta-1} |  \Delta^k_s(u^{(n)}(x,\cdot),t) | ds\\
 &\leq& \frac {C \,t^{-n+\alpha} }{ c^n_{\beta}}\int_t^{\infty} s^{-\beta-1} ds = C_{\alpha, \beta}  \,t^{-n+\alpha-\beta}.
\end{eqnarray*}

Therefore,
$$ \| \frac{\partial^n(P_t  D^\beta f)}{\partial t^n}\|_{\infty}  \leq C t^{-n+\alpha-\beta},$$
and since $\alpha -\beta < n$ by Proposition \ref{Prop1} this  implies $D^\beta f\in Lip_{\alpha-\beta}(\gamma).$\\

\item [ii)] Take $f \in Lip_{\alpha}(\gamma)$, then $f\in L^{\infty}$ and $\|\frac{\partial^n u(\cdot,t)}{\partial
t^n}\|_{\infty}\leq A_{\alpha}(f) t^{-n+\alpha}$. Remember that $\beta < \alpha$,  $k$ is  the smallest integer bigger than $ \beta$ and let $n$ be the smallest integer bigger than $\alpha$.  Observe that from (\ref{Poissboundn}) and (\ref{Lipcond2}), we get,

\begin{eqnarray*}
|  {\mathcal D}^\beta f(x)| &\leq& \frac{1}{c^k_{\beta}}\int_0^{\infty} s^{-\beta-1} |( e^{-s}P_s -I )^k f(x)| \, ds\\
&=&\frac {1}{ c^k_{\beta}}\int_0^1 s^{-\beta-1} |( e^{-s} P_s -I )^k f(x)| ds+ \frac {1}{ c^k_{\beta}}\int_1^{\infty} s^{-\beta-1} |( e^{-s}P_s -I )^k f(x)| ds\\
&=& (I) + (II).
\end{eqnarray*}
Now, let us assume $k<n$, if $k=n$ the argument is straightforward. Let $0 < \varepsilon< 1/2$ such that $\beta+\varepsilon <k$, by Proposition \ref{Prop2} $Lip_{\alpha}(\gamma) \subset Lip_{\beta+\varepsilon}(\gamma)  \subset Lip_{j-\varepsilon}(\gamma), j=1, 2, \cdots, k-1$. Then,  using the identity
$$ (e^{-s} P_s   - I)^k f(x) = \sum_{j=0}^k {k \choose j} e^{-js} (P_s-I)^j (e^{-s}-1)^{k-j} f(x)$$
 and  (\ref{Lipcond2}), we get 
\begin{eqnarray*}
(I) &\leq&\frac {1}{ c^k_{\beta}}\int_0^1 s^{-\beta-1} \sum_{j=0}^k {k \choose j} e^{-js} |(e^{-s}-1)^{k-j}| |(P_s-I)^j f(x)|ds  \\
&\leq&\frac {1}{ c^k_{\beta}} \sum_{j=0}^k {k \choose j} \int_0^1 s^{-\beta-1}e^{-js} |(e^{-s}-1)^{k-j} |
\|(P_s-I)^j f\|_{\infty}\, ds  \\
&\leq&\frac {1}{ c^k_{\beta}}\sum_{j=0}^k {k \choose j} \int_0^1 s^{-\beta-1}  s^{k-j}
\|(P_s-I)^j f\|_{\infty}\, ds  \\
&\leq&  \frac { 1 }{ c^k_{\beta}}  \int_0^1 s^{k-\beta-1}  \, ds\, \|f\|_{\infty}+\sum_{j=1}^{k-1} {k \choose j}\frac {A_{j - \varepsilon}(f)}{ c^k_{\beta}} \int_0^1 s^{k-\beta-1}  s^{k-j}
 s^{j-\varepsilon} \, ds  \\
&& \quad \quad \quad  \quad \quad \quad  \quad \quad \quad \quad   \quad \quad \quad \quad  \quad \quad \quad+\frac {A_{\beta+\varepsilon}(f) }{ c^k_{\beta}} \int_0^1 s^{\beta+\varepsilon-\beta-1} ds\\
&=&  \frac { C }{ c^k_{\beta}(k -\beta)}  \|f\|_{\infty}+  \sum_{j=1}^{k-1} {k \choose j} \frac {A_{j - \varepsilon}(f)}{ c^k_{\beta}(k-\beta-\varepsilon)}+\frac {A_{\beta+\varepsilon}(f) }{ c^k_{\beta}\varepsilon} \\
& \leq& C_{ \beta} \| f \|_{Lip_{\alpha}(\gamma)},
\end{eqnarray*}
and
 \begin{eqnarray*}
(II) &\leq&  \frac {1}{ c_{\beta}}\int_1^{\infty} s^{-\beta-1}[ \sum_{j=0}^k {k \choose j} e^{-(k-j)s}\| P_{(k-j)s} f\|_\infty]ds\\
&\leq & \frac { \|f\|_{\infty}  }{ c^k_{\beta}}\int_1^{\infty} s^{-\beta-1}(1+e^{-s})^k \,ds \leq \frac {2^k \|f\|_{\infty}   }{\beta c_{\beta}}
\leq C_{ \beta} \| f \|_{Lip_{\alpha}(\gamma)}.
\end{eqnarray*}
Thus $D^\beta f \in L^{\infty}$. \\

Now we want to verify the Lipchitz condition.  By Observation \ref{steinobs} it is enough to consider the case $0<t<1$. 
Observe that by the semigroup property,
\begin{eqnarray*}
P_t [( e^{-s}P_s -I )^k f(x) ]  &=&  P_t ( \sum_{j=0}^k {k \choose j} (-1)^je^{-(k-j)s} P_{(k-j)s} f(x) )\\
&=&    \sum_{j=0}^k {k \choose j} (-1)^je^{-(k-j)s} P_{t+(k-j)s} f(x) \\
&=&    \sum_{j=0}^k {k \choose j} (-1)^je^{-(k-j)s} u(x,t+(k-j)s),
\end{eqnarray*}
then 
\begin{eqnarray*}
\frac{\partial^n}{\partial t^n}[ P_t [( e^{-s}P_s -I )^k f(x) ] =  \sum_{j=0}^k {k \choose j} (-1)^je^{-(k-j)s} u^{(n)}(x,t+(k-j)s).
\end{eqnarray*}
Therefore, using  again representation (\ref{kRder}),
\begin{eqnarray*}
\frac{\partial^n(P_t  {\mathcal D}^\beta f(x))}{\partial t^n} &=&  \frac {1}{ c^k_{\beta}}\frac{\partial^n}{\partial t^n} [\int_0^{\infty} s^{-\beta-1} P_t [( P_s -I )^k f(x) ] ds]\\
&=&   \frac {1}{ c^k_{\beta}}\frac{\partial^n}{\partial t^n} [\int_0^{t} s^{-\beta-1} P_t [( P_s -I )^k f(x) ] ds]\\
&&  \quad + \frac {1}{ c^k_{\beta}}\frac{\partial^n}{\partial t^n} [\int_t^{\infty} s^{-\beta-1} P_t [( P_s -I )^k f(x) ] ds]\\
&=&    \frac {1}{ c^k_{\beta}}\int_0^{t} s^{-\beta-1} \sum_{j=0}^k {k \choose j} (-1)^je^{-(k-j)s}  u^{(n)}(x,t+(k-j)s)  ds\\
&&  \quad + \frac {1}{ c^k_{\beta}} \int_t^{\infty} s^{-\beta-1} \sum_{j=0}^k {k \choose j} (-1)^je^{-(k-j)s}  u^{(n)}(x,t+(k-j)s)  ds\\
&=& (III) + (IV).
\end{eqnarray*}

Using  (\ref{difder2}), we have 
\begin{eqnarray*}
 |(III)|    &=&   \frac {e^t}{ c^k_{\beta}}|\int_0^{t} s^{-\beta-1} \sum_{j=0}^k {k \choose j} (-1)^je^{-t-(k-j)s}  u^{(n)}(x,t+(k-j)s)  ds|\\
 &=& \frac {e^t}{ c^k_{\beta}}\int_0^{t} s^{-\beta-1} |  \Delta^k_s(e^{-\cdot}\, u^{(n)}(x,\cdot),t) | ds.\\
\end{eqnarray*} 
Let us take $f(t) =e^{-t}\, u^{(n)}(x,t)$, by (\ref{poisson2}) we know that for any $k >0$
 $$ | \frac{\partial^k( u^{(n)}(x,t) |}{\partial t^k}|  \leq C t^{-(n+k)+\alpha} = C t^{-n+(\alpha-k)} ,$$
 and then using Leibnitz formula, and the fact that $0<t<1$
 \begin{eqnarray*}
 |\frac{\partial^k[ e^{-t}u^{(n)}(x,t))]}{\partial t^k}| &=& | e^{-t} \sum_{j=0}^k {k \choose j} (-1)^{(j)} u^{n+(k-j)}(x,t)| \\
&\leq&  e^{-t} \sum_{j=0}^k {k \choose j}|u^{n+(n-j)}(x,t)| \\
&\leq&  C e^{-t} \sum_{j=0}^k {k \choose j} t^{-(n+(k-j))+\alpha} \\
&=& C e^{-t} t^{-n+\alpha} \sum_{j=0}^k {k \choose j} t^{-(k-j))} \\
&\leq& C e^{-t} t^{-n+\alpha} 2^k t^{-k} = C e^{-t} t^{-(n+k)+\alpha}
  \end{eqnarray*}
  Then with a small variation  of the argument of the proof of (\ref{difder1}) we get
  \begin{eqnarray*}
|  \Delta^k_s(e^{-\cdot}u^{(n)}(x,\cdot),t) | &\leq&  C t^{-(n+k)+\alpha} e^{-t} s^k,
\end{eqnarray*}
and therefore
 $$|(III)| \leq  \frac { C\, t^{-(n+k)+\alpha}}{ c^k_{\beta}}\int_0^{t} s^{-\beta+k-1} ds = C_{\alpha, \beta}  \, t^{-k+\alpha-n} t^{-\beta+k} = C_{\alpha, \beta}  \, t^{-n+\alpha-\beta}.$$

On the other hand, 
\begin{eqnarray*}
 |(IV)| &\leq&  \frac {1}{ c^k_{\beta}} \int_t^{\infty} s^{-\beta-1} \sum_{j=0}^k {k \choose j} e^{-(k-j)s}  |u^{(n)}(x,t+(k-j)s)|  ds\\
 &\leq&  \frac {1}{ c^k_{\beta}} \int_t^{\infty} s^{-\beta-1} 2^k  (t+(k-j)s)^{-n+\alpha} ds\\
 &\leq& C \,t^{-n+\alpha} \int_t^{\infty} s^{-\beta-1} ds=C_{\alpha, \beta}  \,t^{-n+\alpha-\beta}.
\end{eqnarray*}
Therefore, we have we can conclude that

$$ \| \frac{\partial}{\partial t}(P_t  {\mathcal D}^\beta f)\|_{\infty}  \leq C_{\alpha, \beta} t^{-n+\alpha-\beta}, $$
and again, since $\alpha -\beta < n$ by Proposition \ref{Prop1} this  implies  $ {\mathcal D}^\beta f\in Lip_{\alpha-\beta}(\gamma).$ \ep \\

\end{enumerate}

\begin{obs}
 Let us observe that the arguments given in the proofs of Theorem 2.1 and 2.2 are valid in the classical case taking the Poisson integral, and therefore they are alternative proofs of the ones given in \cite{se70}.
\end{obs}
\begin{obs}
 Moreover, if instead of considering the {\em Ornstein-Uhlenbeck operator} (\ref{OUop}) and the {\em Poisson-Hermite semigroup} (\ref{PoissonH}) we consider the {\em  Laguerre differential  operator} in $\mathbb R^d_{+}$.
\begin{equation}
\mathcal{L}^{\alpha} = \sum^{d}_{i=1} \bigg[ x_i \frac{\partial^2}{\partial x^2_i}
+ (\alpha_i + 1 - x_i )  \frac{\partial}{\partial x_i} \bigg],
\end{equation}
and the corresponding {\em Poisson-Laguerre semigroup}, or if we consider the {\em Jacobi differential  operator} in $(-1,1)^d$,
\begin{equation}
\mathcal{L}^{\alpha,\beta} = - \sum^{d}_{i=1} \bigg[ (1-x_i^2)\frac{\partial^2}{\partial x_i^2}
+  (\beta_i -\alpha_i-\left(\alpha_i +\beta_i +2\right)x_i) \frac{\partial}{\partial x_i} \bigg],
\end{equation}
and the corresponding {\em Poisson-Jacobi semigroup} (for details we refer to \cite{ur2}), the arguments are completely analogous.  That is to say, we can defined in analogous manner {\em Laguerre-Lipschitz spaces} and {\em Jacobi-Lipschitz spaces}  and prove that the corresponding notions of Fractional Integrals and Fractional Derivatives, see \cite{GLLNU}, \cite{balur}  behave similarly.
In order to see this it is more convenient to use the representation (\ref{PoissonH}) of $P_t$ in terms of the one-sided stable measure $\mu^{(1/2)}_t(ds)$ and the write Lemma \ref{L1estimate} in terms of it, see \cite{PinUr}.
\end{obs}

{\bf Acknowldegment}
We want to thanks Ebner Pineda for his corrections to various version of the paper as well as some interesting remarks.

\section{Appendix}

Proof of Lemma \ref{lemaforwdiff}

\begin{enumerate}
\item[i)] Let us prove the first equality, the second one is totally analogous,
\begin{eqnarray*}
\Delta_s^{k-1} (\Delta_s(f,\cdot),t)&=& \sum_{j=0}^{k-1} {k-1 \choose j} (-1)^j \Delta_s(f, t+(k-1-j)s)\\
&=& \sum_{j=0}^{k-1} {k-1 \choose j} (-1)^j f( t+(k-j)s)\\
&& \quad \quad \quad  - \sum_{j=0}^{k-1} {k-1 \choose j} (-1)^j f(t+(k-1-j)s)
\end{eqnarray*}
\begin{eqnarray*}
\quad \quad \quad &&\\
&=& f(t+ks) + \sum_{j=1}^{k-1} {k-1 \choose j} (-1)^j f( t+(k-j)s)\\
&& \quad + \sum_{j=0}^{k-2} {k-1 \choose j} (-1)^{(j+1)} f(t+(k-(j+1))s) +(-1)^k f(t)\\
&=& f(t+ks) + \sum_{j=1}^{k-1} {k-1 \choose j} (-1)^j f( t+(k-j)s)\\
&& \quad + \sum_{j=1}^{k-1} {k-1 \choose j+1} (-1)^{j} f(t+(k-j)s) +(-1)^k f(t)\\
&=& f(t+ks) + \sum_{j=1}^{k-1} [{k-1 \choose j} +{k-1 \choose j+1}](-1)^j f( t+(k-j)s)+(-1)^k f(t)\\
&=& f(t+ks) + \sum_{j=1}^{k-1} {k \choose j} (-1)^j f( t+(k-j)s)+(-1)^k f(t) = \Delta^k_s(f,t).
\end{eqnarray*}

\item[ii)] By induction in $k$. For $k=1$, using the Fundamental Theorem
$$ \Delta_s(f,t) = f(t+s) - f(t) =\int_t^{t+s} f'(v) dv.$$
Let us assume that the identity is true for $k-1$,
$$\Delta^{k-1}_s(f,t) = \int_t^{t+s} \int_{v_1}^{v_1+s} \cdots \int_{v_{k-2}}^{v_{k-2}+s} f^{(k-1)}(v_{k-1}) \,  dv_{k-1} \cdots dv_2 dv_1,$$
and let us prove it for $k$. Using i) and the Fundamental Theorem, we get,  after performing $k-1$ change of variables,

\begin{eqnarray*}
\Delta^k_s(f,t) &=&  \Delta_s(\Delta_s^{k-1}(f,\cdot),t) = \Delta_s^{k-1}(f,t+s) -\Delta_s^{k-1}(f,t)\\
&=& \int_{t+s}^{t+2s} \int_{v_1}^{v_1+s} \cdots \int_{v_{k-2}}^{v_{k-2}+s} f^{(k-1)}(v_{k-1}) \,  dv_{k-1} \cdots dv_2 dv_1 \\
&& \quad \quad  - \int_t^{t+s} \int_{v_1}^{v_1+s}\cdots \int_{v_{k-2}}^{v_{k-2}+s} 
f^{(k-1)}(v_{k-1}) \,  dv_{k-1} \cdots dv_2 dv_1\\
&=& \int_{t}^{t+s} \int_{v_1}^{v_1+s} \cdots \int_{v_{k-2}+s}^{v_{k-2}+2s} f^{(k-1)}(v_{k-1}) \,  dv_{k-1} \cdots dv_2 dv_1 \\
&& \quad \quad  - \int_t^{t+s} \int_{v_1}^{v_1+s} \cdots \int_{v_{k-2}}^{v_{k-2}+s} 
f^{(k-1)}(v_{k-1}) \,  dv_{k-1} \cdots dv_2 dv_1\\
&=& \int_{t}^{t+s} \int_{v_1}^{v_1+s} \cdots [\int_{v_{k-2}+s}^{v_{k-2}+2s} f^{(k-1)}(v_{k-1}) \,  dv_{k-1} \\
&& \quad \quad  \quad \quad  \quad \quad  \quad \quad -\int_{v_{k-2}}^{v_{k-2}+s} 
f^{(k-1)}(v_{k-1}) \,  dv_{k-1}] \cdots dv_2 dv_1\cdots dv_2 dv_1 \\
&=& \int_{t}^{t+s} \int_{v_1}^{v_1+s} \cdots \int_{v_{k-2}}^{v_{k-2}+s} [f^{(k-1)}(v_{k-1}+s) -f^{(k-1)}(v_{k-1}) ]\,  dv_{k-1} \cdots dv_2 dv_1 \\
&=& \int_{t}^{t+s} \int_{v_1}^{v_1+s}  \cdots \int_{v_{k-2}}^{v_{k-2}+s}\int_{v_{k-1}}^{v_{k-1}+s} f^{(k)}(v_{k}) \,  dv_{k} dv_{k-1} \cdots dv_2 dv_1.
\end{eqnarray*}
\item[iii)] Let us  prove (\ref{difder}), 
\begin{eqnarray*}
\frac{\partial }{\partial s}(\Delta_s^k (f,t) )&=& D_s(\sum_{j=0}^k {k \choose j} (-1)^j f(t+(k-j)s))\\
& =&  \sum_{j=0}^k {k \choose j} (-1)^j \frac{\partial }{\partial s}(f(t+(k-j)s)) \\
& =&  \sum_{j=0}^{k-1} {k \choose j} (-1)^j (k-j) f'(t+(k-j)s) \\
& =& k \sum_{j=0}^{k-1} {k-1 \choose j} (-1)^jf'((t+s)+(k-1-j)s) \\
&=&k \,\Delta_s^{k-1} (f',t+s).
\end{eqnarray*}
Now, let us prove  (\ref{difder2})
\begin{eqnarray*}
\frac{\partial^j }{\partial t^j}(\Delta_s^k (f,t) )&=& \frac{\partial^j }{\partial t^j}(\sum_{j=0}^k {k \choose j} (-1)^j f(t+(k-j)s))\\
& =&  \sum_{j=0}^k {k \choose j} (-1)^j \frac{\partial^j }{\partial t^j}(f(t+(k-j)s)) \\
& =&  \sum_{j=0}^k {k \choose j} (-1)^j f^{(j)}(t+(k-j)s)) \\
&=& \,\Delta_s^{k} (f^{(j)},t).
\end{eqnarray*}
\end{enumerate}
\ep


\begin{thebibliography}{99}
\bibitem{balur}
Balderrama,C., Urbina, W. \emph{ Fractional Integration and Fractional Differen-tiation for d-dimensional Jacobi Expansions.}  Contemporary Math. AMS \#471 (2008) 1-14.

\bibitem{calderon}
 Calder\'on, A. P. \emph{ Intermediate spaces and interpolation, the complex method} Studia Math. 24 (1964), 113Ð190.

\bibitem{fort}
Fort, T. \emph{Finite Differences and difference Equations in the Real Domain} Oxford Claderon Press (1948).
 
 \bibitem{forscotur}
Forzani, L., Scotto, R, and Urbina, W{.} \emph{Riesz and Bessel
Potentials, the $g_k$ functions and an Area function, for  the
Gaussian measure $\gamma_d$.} Revista de la Uni\'on Matem\'atica
 Argentina (UMA), vol 42 (2000), no.1,17--37.

\bibitem{fm91}
  Frazier M., Jawerth B., Weiss G. \emph{Littlewood Paley Theory and the Study of Functions Spaces.} CBMS-Conference Lecture Notes 79.
 Amer. Math. Soc. Providence RI , 1991.
 
 \bibitem{GLLNU}
 Graczyk, P, Loeb, J.J., L\'opez, I., Nowak, A.  \& Urbina, W \emph{ Higher order Riesz transforms,  Fractional differentiation and Sobolev spaces for Laguerre expansions.}  J. Math. Pures Appl. (9). 84 (2005), no. 3, 375-405
 
 \bibitem{gs96}
 Gatto, A. E, Segovia, C., V\'{a}gi, S. \emph{On Fractional Differentiation and Integration on Spaces of Homogeneous Type}, Rev. Mat. 
Iberoamericana 12 (1996) 111-145.
  
\bibitem{gatto08}
 Gatto, A. \emph{Boundedness on inhomogeneous Lipschitz spaces of fractional integrals, singular integrals and hypersingular integrals associated to non-doubling measures} Collect. Math. 60, 1 (2009), 101Ð114.
 
 \bibitem{Heideman}
 Heideman, N. J. H. \emph{Duality and fractional integration on Lipschitz spaces} Studia Math.  50 (1974), 65Ð85.

\bibitem{lour}
L\'{o}pez I{.} and  Urbina, W{.} \emph{Fractional Differentiation
for the Gaussian Measure and Applications. } Bull. Sciences Math, 2004,
{\bf 128}, 587--603.

\bibitem{PinUr}
Pineda, E., Urbina, W. \emph{Some results on gaussian Besov-Lipschitz spaces and Gaussian Triebel-Lizorkin spaces.} Journal of Approximation Theory (2008), doi:10.1016/j.jat.2008.11.010 In Press.

 \bibitem{sam}
  Samko S., Kilbas, A \& Marichev, O.  \emph{Fractional integrals and derivatives: theory and applications.} Gordon and Breach Science Publishers, Philadelphia, 1992.

\bibitem{sp97}
  Sj\"{o}gren P. \emph{Operators associated with the Hermite semigroup- a survey.} J. Fourier Anal. Appl. {\bf 3} (1997), Special Issue, 813--823.
  
\bibitem{se70}
  Stein E. \emph{Singular integrals and differentiability properties of functions} Princeton Univ. Press. Princeton, New Jersey, 1970.
  
  \bibitem{st2}
 Stein, E. M. {\em Topics in Harmonic Analysis related to the Littlewood-Paley
 Theo-ry.} Princeton Univ. Press. Princeton (1970).
  
   \bibitem{trie}
 Triebel, H  \emph{Theory of function spaces II.} Birkh\"auser Verlag, Basel, 1992.
 
  \bibitem{ur2}
 Urbina, W.  {\em Operators Semigroups associated to Classical Orthogonal Polynomials  and  Functional Inequalities.} Lecture Notes of the French Mathematical Society (SMF). 2008.
\end{thebibliography}
\end{document}